\documentclass[11pt,intlimits,reqno]{amsart}
\usepackage{graphicx}
\usepackage{pb-diagram}
\usepackage{pstricks}
\usepackage[all]{xy}
\usepackage{a4wide}
\usepackage{hyperref}
\usepackage{amsfonts,amssymb}
\usepackage{mathrsfs}
\usepackage{txfonts}


\title[Survival of infinitely many critical points]{Survival of infinitely many critical points \\for the Rabinowitz action functional}

\author[J. Kang]{Jungsoo Kang}
\address{Department of Mathematical Sciences, Seoul National University,
Kwanakgu Shinrim, San56-1 Seoul, South Korea}
\email{hoho159@snu.ac.kr}

\newtheorem{neu}{}[section]

\newtheorem{Prop}[neu]{Proposition}
\theoremstyle{definition}

\newtheorem{Rmk}[neu]{Remark}
\newtheorem{Def}[neu]{Definition}

\theoremstyle{remark}

\theoremstyle{definition}

\newcommand{\p}{\partial}
\newcommand{\om}{\omega}

\newcommand{\pf}{\longrightarrow}

\newcommand{\Z}{{\mathbb{Z}}}
\newcommand{\R}{{\mathbb{R}}}

\renewcommand{\H}{\mathrm{H}}

\newcommand{\im}{\mathrm{im}}
\newcommand{\CF}{\mathrm{CF}}
\newcommand{\HF}{\mathrm{HF}}
\newcommand{\RFH}{\mathrm{RFH}}

\newcommand{\Crit}{{\rm Crit}}

\newcommand{\Ham}{\mathrm{Ham}}
\newcommand{\Supp}{\mathrm{Supp}}

\renewcommand{\AA}{\mathcal{A}}
\newcommand{\x}{\times}

\newcommand{\beq}{\begin{equation}}
\newcommand{\beqn}{\begin{equation}\nonumber}
\newcommand{\eeq}{\end{equation}}

\newcommand{\bea}{\begin{equation}\begin{aligned}}
\newcommand{\bean}{\begin{equation}\begin{aligned}\nonumber}
\newcommand{\eea}{\end{aligned}\end{equation}}

\numberwithin{equation}{section}

\begin{document}
%%%%%%%%%%%%%%%%%%%%%%%%%%%%%%%%%%%%%%%%%%%%

\begin{abstract}
In this paper, we show that if Rabinowitz Floer homology has infinite dimension, there exist infinitely many critical points of a Rabinowitz action functional even though it could be non-Morse. This result is proved by examining filtered Rabinowitz Floer homology.
\end{abstract}

\maketitle

\section{Introduction}
Many recent studies have focused on the coisotropic intersection problem. In particular, there has been growing interest in the leafwise intersection property described below.\\[-1ex]

\noindent\textbf{Leafwise intersection.}
Let $(M,\om)$ be a $2n$ dimensional symplectic manifold and $\Sigma$ be a coisotropic submanifold of codimension $k$. Then the symplectic structure $\om$ determines a symplectic orthogonal bundle $T\Sigma^\om$, which is a subbundle of $T\Sigma$ by the definition of coisotropic:
\beqn
T\Sigma^\om:=\{(x,\xi)\in T\Sigma\,|\,\om_x(\xi,\zeta)=0 \textrm{ for all } \zeta\in T_x\Sigma\}
\eeq
Since $\om$ is closed, $T\Sigma^\om$ is integrable, thus $\Sigma$ is foliated by the leaves; we denote by $L_x$ the leaf through $x\in\Sigma$. We call $x\in\Sigma$ a {\em leafwise intersection} of $\phi\in\Ham_c(M,\om)$ if $x\in L_x\cap\phi(L_x)$. $\Ham_c(M,\om)$ is defined below in Conventions and Notations.

In this article, we focus on the case that $\Sigma$ is a restricted contact type hypersurface with contact form $\lambda$, i.e. $d\lambda=\om$ and $\lambda\wedge\om^{n-1}|_\Sigma\neq0$, which bounds a compact region in $M$. The argument developed in this article would continue to hold in general contact case with minor modification (see \cite{Ka1} for the generalized Rabinowitz Floer theory). In that case, $T\Sigma^\om$ is nothing but the characteristic line bundle spanned by the Reeb vector field $R$. A leafwise intersection $x\in\Sigma$ is called {\em a periodic leafwise intersection} if the leaf $L_x$ through $x$ is a closed Reeb orbit.\\[-1ex]

The problem of finding leafwise intersection points was initiated by Moser \cite{Mo} and pursued further in \cite{Ba,Dr,EH}. Recently there are three different techniques to this problem; \cite{Zi} approached the leafwise intersection problem with Lagrangian Floer homology and \cite{Gi,Gu} took advantage of Symplectic Floer homology. Another effective way to investigate leafwise intersections is Rabinowitz Floer homology developed by Cieliebak-Frauenfelder \cite{CF}. Albers-Frauenfelder \cite{AF1} observed that critical points of a perturbed Rabinowitz action functional give rise to leafwise intersections; \cite{AF1,AF2,AF3,AF4,AMo,AMc,Ka,Ka1,Ka2,Me} obtained many results on the leafwise intersection problem using Rabinowitz Floer homology. Rabinowitz Floer homology is the Floer homology of the Rabinowitz action functional $\AA^H_F$. Thus, the dimension of Rabinowitz Floer homology gives a lower bound on the number of leafwise intersections as in the Morse inequalities in case the Rabinowitz action functional is Morse. Interestingly enough, Rabinowitz Floer homology has infinite dimension in some cases \cite{CFO,Me,AMc,Ka2}. In such cases, we know that for $F\in\Ham_c(M,\om)$ which makes $\AA^H_F$ Morse, $\AA^H_F$ has infinitely many critical points. It turned out that the Rabinowitz action functional is Morse for a generic $F\in\Ham_c(M,\om)$ \cite{AF1}; therefore a generic perturbation has either infinitely many leafwise intersections or a periodic leafwise intersection (see Proposition \ref{prop:critical point answers question}); furthermore, \cite{AF2} showed that generically there is no periodic leafwise intersection points.\\[-1ex]

\noindent\textbf{Main Theorem.} {\em If Rabinowitz Floer homology has infinite dimension, then for $\phi_F\in\Ham_c(M,\om)$ there exists infinitely many leafwise intersections or a periodic leafwise intersection.}\\[-1ex]

It is noteworthy that Main Theorem does not assume any kind of non-degeneracy for $\phi_F$; that is, the perturbed Rabinowitz action functional $\AA^H_F$ is not necessarily Morse. When it is Morse, the theorem follows immediately from the Morse inequalities.
\begin{Rmk}
An analogous result continues to hold without doubt in any other Morse or Floer homology whenever the action functional is a Morse (resp. Morse-Bott) and it has only finitely many critical points (resp. finitely many compact critical components) in a compact action interval.
\end{Rmk}
\begin{Rmk}
The Main theorem subsumes a result in \cite{AF3} which showed the above result for a (unit) cotangent bundle by means of the spectral invariants in Rabinowitz Floer homology.
\end{Rmk}

\noindent\textbf{Convention and Notations.}
\begin{itemize}
\item The {\em Reeb vector field} $R$ is characterized by $\lambda(R)=1$ and $i_R d\lambda=0$.
\item The {\em Hamiltonian vector field} $X_F$ associated to a Hamiltonian function $F\in C^\infty(S^1\x M)$ is defined implicitly by $i_{X_F}\om=dF$.
\item $\phi_F$ is the time one map of the flow of $X_F$ and called a {\em Hamiltonian diffeomorphism}.
\item We denote by $\Ham_c(M,\om)$ the group of Hamiltonian diffeomorphisms generated by compactly supported Hamiltonian functions.
\end{itemize}

\subsection{Filtered Rabinowitz Floer homology}
Since $\Sigma$ is a contact hypersurface, there exists a Liouville vector field $Y$ such that $L_Y\om=\om$ and $Y\pitchfork\Sigma$; we denote by $\phi_Y^t$ the flow of $Y$ and fix $\delta>0$ such that $\phi_Y^t|_\Sigma$ is defined for $|t|<\delta$. Since $\Sigma$ bounds a compact region in $M$, we are able to define a Hamiltonian function $G\in C^\infty(M)$ so that
\begin{enumerate}
\item $G^{-1}(0)=\Sigma$ is a regular level set;
\item $G(\phi_Y^t(x))=t$ for all $x\in\Sigma$ and $|t|<\delta$;
\item $dG$ has compact support.
\end{enumerate}

\begin{Def}
Given time-dependent Hamiltonian functions $H\in C^\infty(S^1\x M)$ and $F\in C_c^\infty(S^1\x M)$, a pair $(H,F)$ is called a {\em Moser pair} if it satisfies
\begin{enumerate}
\item $H$ is a weakly time-dependent Hamiltonian function. That is, $H$ is of the form $H(t,x)=\chi(t)G(x)$ for $G\in C^\infty(M)$ defined above and $\chi:S^1\to S^1$ with $\int_0^1\chi dt=1$ and $\Supp\chi\subset(\frac{1}{2},1)$;
\item their time supports are disjoint, i.e.
\beqn
H(t,\cdot)=0 \quad\textrm{for}\,\,\, \forall t\in\big[0,\frac{1}{2}\big] \quad\textrm{and}\quad F(t,\cdot)=0 \quad\textrm{for}\,\,\, \forall t\in\big[\frac{1}{2},1\big].
\eeq
\end{enumerate}
\end{Def}
\begin{Rmk}
We can easily check that every element in $\Ham_c(M,\om)$ is generated by a compactly supported Hamiltonian function with time support on $\big[0,\frac{1}{2}\big]$.
\end{Rmk}

We define a perturbed Rabinowitz Floer functional $\AA^H_F:C^\infty(S^1,M)\x\R\pf\R$ for a Moser pair $(H,F)$ as follows:
\beqn
\AA^H_F(v,\eta)=-\int_0^1v^*\lambda-\eta\int_0^1H(t,v)dt-\int_0^1F(t,v)dt.
\eeq
Albers-Frauenfelder observed that a critical point of $\AA^H_F$ gives rise to a leafwise intersection.
\begin{Prop}\label{prop:critical point answers question}
\cite{AF1} Let $(v,\eta)\in\Crit\AA_{F}^{H}$. Then $x:=v(0)\in\Sigma$ and satisfies $\phi_F(x)\in L_x$. Thus, $x$ is a leafwise intersection point. Moreover, the map
\beqn
\Crit\AA^{H}_F\pf\big\{\textrm{leafwise intersections}\big\}
\eeq
is injective unless there exists a periodic leafwise intersection.
\end{Prop}

\begin{Def}
We define the {\em action spectrum} for a Moser pair $(H,F)$ by
\beqn
\mathrm{Spec} (H,F):=\big\{\AA^H_F(v,\eta)\in\R\,\big|\,(v,\eta)\in\Crit(\AA^H_F)\big\}.
\eeq
\end{Def}

Now, we roughly describe a filtered Rabinowitz Fleor homology for a Morse action functional $\AA^H_F$ and $a<b<c\notin\mathrm{Spec}(H,F)$ (see \cite{Ka1} for a rigorous definition); it is proved in \cite{AF1} that $\AA^H_F$ is Morse for a generic $F\in C_c^\infty(S^1\x M)$ with time support in $[0,1/2]$. Let us set the $\Z/2$-module by
\beqn
\CF^{(a,b)}(\AA^H_F):=\Crit^{(a,b)}(\AA^H_F)\otimes\Z/2
\eeq
where
\beqn
\Crit^{(a,b)}(\AA^H_F):=\big\{(v,\eta)\in\Crit\AA^H_F\,\big|\,\AA^H_F(v,\eta)\in(a,b)\big\}.
\eeq
Then it becomes a complex with the boundary operator $\p$ defined by counting solutions of a nonlinear elliptic PDE. Then the filtered Rabinowitz Floer homology is defined by
\beqn
\HF^{(a,b)}(\AA^H_F)=\H\big(\CF^{(a,b)}(\AA^H_F),\p\big).
\eeq
In the case $a=-\infty$ and $b=\infty$, we get a full Rabinowitz Floer homology $\HF(\AA^H_F)$ and it is invariant under the choice of a perturbation $F\in C_c^\infty(S^1\x M)$ by the standard argument in the Floer theory. Thus we denote by
\beqn
\RFH(\Sigma,M):=\HF(\AA^H)\cong\HF(\AA^H_F)
\eeq
where $\AA^H=\AA^H_0$. The Rabinowitz action functional $\AA^H$ is Morse-Bott for a generic contact hypersurface; in that case Rabinowitz Floer homology $\RFH(\Sigma,M)$ is well-defined. Otherwise, we define $\RFH(\Sigma,M):=\HF(\AA^H_f)$ for $\AA^H_f$ being Morse. Moreover we have canonical homomorphisms
\beqn
i_a^{b,c}:\CF^{(a,b)}(\AA^H_F)\pf\CF^{(a,c)}(\AA^H_F),\quad \pi_{a,b}^c:\CF^{(a,c)}(\AA^H_F)\pf\CF^{(b,c)}(\AA^H_F)
\eeq
inclusions and projections respectively. They fit into a short exact sequence
\beqn
0\pf\CF^{(a,b)}(\AA^H_F)\stackrel{i_a^{b,c}}{\pf}\CF^{(a,c)}(\AA^H_F)\stackrel{\pi^c_{a,b}}{\pf}\CF^{(b,c)}(\AA^H_F)\pf0,
\eeq
and we obtain a long exact sequence
\beqn
\cdots\stackrel{\pi_{a,b*}^c}{\pf}\HF^{(b,c)}(\AA^H_F)\stackrel{\delta_*}{\pf}\HF^{(a,b)}(\AA^H_F)\stackrel{i_{a*}^{b,c}}{\pf}\HF^{(a,c)}(\AA^H_F)\stackrel{\pi_{a,b*}^c}
{\pf}\HF^{(b,c)}(\AA^H_F)\stackrel{\delta_*}{\pf}\cdots.
\eeq

Next, we briefly recall the definition of Hofer norm.
\begin{Def}\label{def:Ham}
Let $F\in C^\infty_c(S^1\x M,\R)$ be a compactly supported time-dependent Hamiltonian function. We set
$$
||F||_+:=\int_0^1\max_{x\in M} F(t,x) dt\qquad||F||_-:=-\int_0^1\min_{x\in M} F(t,x) dt=||-F||_+
$$
and
$$
||F||:=||F||_++||F||_-.\;
$$
For $\phi\in \Ham_c(M,\om)$ the Hofer norm is
$$
||\phi||:=\inf\{||F||\mid \phi=\phi_F,\,\, F\in C^\infty_c(S^1\x M,\mathbb{R})\}.\;
$$

\end{Def}

The following proposition is well-known in the standard Floer theory and easily shown by energy estimates (see e.g. \cite{Ka1}).
\begin{Prop}\label{prop:continuation}
There exist continuation homomorphisms for the filtered case:
\bean
\Phi&:\HF^{(a,b)}(\AA^H_f)\pf\HF^{(a+||F-f||_-,b+||F-f||_-)}(\AA^H_F),\\
\Psi&:\HF^{(a,b)}(\AA^H_F)\pf\HF^{(a+||F-f||_+,b+||F-f||_+)}(\AA^H_f).
\eea
\end{Prop}

\section{Proof of Main Theorem}
\begin{proof}
Assume on the contrary that $\RFH(\Sigma,M)=\HF(\AA^H_f)$ is infinite dimensional and $\AA^H_F$ for $F\in C^\infty_c(S^1\x M)$ has finitely many critical points. Then there exist $A,B\in\R$ such that
$$
\mathrm{Spec}(H,F)\subset[A,B].
$$
Since $F$ and $f$ are compactly supported, $||F-f||_-$ has finite value, so we can pick $b\in\R$ so that $b>B-||F-f||_-$. Then $\CF^{(b+||F-f||_-,\infty)}(\AA^H_F)$ is empty and thus its homology $\HF^{(b+||F-f||_-,\infty)}(\AA^H_F)$ is well-defined although $\AA^H_F$ could have degenerate critical points with action outside $(b+||F-f||_-,\infty)$. Due to Proposition \ref{prop:continuation} together with the ``homotopy of homotopies'' argument, we have the following commutative diagram \eqref{eq:diagram}.
\bea\label{eq:diagram}
\xymatrix{ \HF^{(b,\infty)}(\AA^H_f) \ar[dd]^{\big(\pi^\infty_{b,b+||F-f||}\big)*} \ar[dr]^\Phi & \hspace{8cm}  \\
 &  \HF^{(b+||F-f||_-,\infty)}(\AA^H_F) \ar[dl]^\Psi \\
\HF^{(b+||F-f||,\infty)}(\AA^H_f) & \hspace{8cm} }
\eea
Since $\HF^{(b+||F||_-,\infty)}(\AA^H_F)$ vanishes, accordingly $\big(\pi^\infty_{b,b+||F-f||}\big)*$ also vanishes and we have the following short exact sequence:
\beqn
0\stackrel{\pi_*}{\pf}\HF^{(b+||F-f||,\infty)}(\AA^H_f)\stackrel{\delta_*}{\pf}\HF^{(b,b+||F-f||)}(\AA^H_f)\stackrel{i_*}
{\pf}\HF^{(b,\infty)}(\AA^H_f)\stackrel{\pi_*}{\pf}0
\eeq
where $\pi_*=\big(\pi^\infty_{b,b+||F-f||}\big)_*$, $i_*=\big(i_{b}^{b+||F-f||,\infty}\big)_*$, and $\delta_*$ is a connecting homomorphism. Since $||F||<\infty$, the dimension of $\HF^{(b,b+||F-f||)}(\AA^H_f)$ is finite. Since $i_*$ is surjective, $\HF^{(b,\infty)}(\AA^H_f)$ is finite dimensional.

In the same way, we choose $a\in\R$ so that $a<A-||F-f||_-$ and then the following diagram commutes.
\bean
\xymatrix{ \HF^{(-\infty,a)}(\AA^H_f) \ar[dd]^{\big(i_{-\infty}^{a,a+||F-f||}\big)_*} \ar[dr]^\Phi & \hspace{8cm}  \\
 &  \HF^{(-\infty,a+||F-f||-)}(\AA^H_F) \ar[dl]^\Psi \\
\HF^{(-\infty,a+||F-f||)}(\AA^H_f) & \hspace{8cm} }
\eea
Analogously, $\big(i_{-\infty}^{a,a+||F-f||}\big)_*$ vanishes and thus we have
\beqn
0\stackrel{i_*'}{\pf}\HF^{(-\infty,a+||F-f||)}(\AA^H_f)\stackrel{\pi_*'}{\pf}\HF^{(a,a+||F-f||)}(\AA^H_f)\stackrel{\delta_*}{\pf}\HF^{(-\infty,a)}(\AA^H_f)
\stackrel{i_*'}{\pf}0.
\eeq
where $i_*'=\big(i_{-\infty}^{a,a+||F-f||}\big)_*$, $\pi_*'=\big(\pi^{a+||F-f||}_{-\infty,a}\big)_*$, and $\delta_*$ is a connecting homomorphism.
Therefore $\HF^{(-\infty,a)}(\AA^H_f)$ has finite dimension. Next, we compare the dimensions of homologies and this is contradictory to the our assumption that $\AA^H_F$ has finitely many critical points. Consider a part of the long exact sequence:
\beqn
\cdots\stackrel{\delta_*}{\pf}\HF^{(b,\infty)}(\AA^H_f)\stackrel{i_{b,a*}^\infty}{\pf}\HF^{(a,\infty)}(\AA^H_f)\stackrel{\pi_{a*}^{\infty,b}}{\pf}\HF^{(a,b)}
(\AA^H_f)\stackrel{\delta_*}{\pf}\cdots.
\eeq
Then we know that
\bea\label{eq:dim cal}
\dim\HF^{(a,\infty)}(\AA^H_f)&=\dim\im\,\pi_{a*}^{\infty,b}+\dim\ker\pi_{a*}^{\infty,b}\\
&=\dim\im\,\pi_{a*}^{\infty,b}+\dim\im \,i_{b,a*}^\infty\\
&\leq\dim\HF^{(a,b)}(\AA^H_f)+\dim\HF^{(b,\infty)}(\AA^H_f)<\infty.
\eea
In addition, we have another long exact sequence of the form
\beq\label{eq:les}
\cdots\stackrel{\delta_*}{\pf}\HF^{(-\infty,a)}(\AA^H_f)\stackrel{i_{-\infty*}^{a,\infty}}{\pf}\RFH(\Sigma,M)\stackrel{\pi_{-\infty,a*}^\infty}{\pf}
\HF^{(a,\infty)}(\AA^H_f)\stackrel{\delta_*}{\pf}\cdots.
\eeq
Using the computation \eqref{eq:dim cal} together with \eqref{eq:les}, we derive the following contradiction and it finishes the proof of theorem.
\bean
\infty=\dim\RFH(\Sigma,M)&=\dim\im\,\pi_{-\infty,a*}^\infty+\dim\ker\pi_{-\infty,a*}^\infty\\
&=\dim\im\,\pi_{-\infty,a*}^\infty+\dim\im \,i_{-\infty*}^{a,\infty}\\
&\leq\dim\HF^{(a,\infty)}(\AA^H_f)+\dim\HF^{(-\infty,a)}(\AA^H_f)<\infty.
\eea
\end{proof}

\subsubsection*{Acknowledgments}
I owe a debt of gratitude to Urs Frauenfelder for valuable discussions. This work is partially supported by the Basic Research grant Nr. 2010-0007669 from the Korea government.

\end{document}